\documentclass[12pt]{article}

\usepackage[left=1in,top=1in,right=1in,bottom=1in]{geometry}

\usepackage{amsmath,amssymb,amsthm,amsfonts}
\usepackage{times}
\usepackage{enumerate}
\usepackage{epsfig}
\usepackage{graphicx}
\usepackage{pgf}

\newtheorem{theorem}{Theorem}[section]
\newtheorem*{theorem*}{Theorem}

\newtheorem*{corollary*}{Corollary}

\newtheorem*{conjecture*}{Conjecture}

\theoremstyle{remark}

\newtheorem*{remark}{Remark}

\theoremstyle{definition}
\newtheorem{definition}[theorem]{Definition}

\newtheorem{example}[theorem]{Example}
\numberwithin{figure}{section}

\def\R{\mathbb R}
\def\Z{\mathbb Z}
\def\C{\mathcal{C}}
\def\S{\mathcal{S}}
\def\Re{\mathcal{R}}

\def\E{\mathbb{E}}

\begin{document}

\title{Product-form stationary distributions for deficiency zero
  chemical reaction networks} \author{David F.  Anderson,$^{1}$
  Gheorghe Craciun,$^{1}$ and Thomas G. Kurtz$^{1}$}

\footnotetext[1]{Department of Mathematics, University of Wisconsin,
  Madison, WI, 53706}

\maketitle

\begin{abstract}
  We consider stochastically modeled chemical reaction systems with
  mass-action kinetics and prove that a product-form stationary
  distribution exists for each closed, irreducible subset of the state
  space if an analogous deterministically modeled system with
  mass-action kinetics admits a complex balanced equilibrium.
  Feinberg's deficiency zero theorem then implies that such a
  distribution exists so long as the corresponding chemical network is
  weakly reversible and has a deficiency of zero.  The main parameter
  of the stationary distribution for the stochastically modeled system
  is a complex balanced equilibrium value for the corresponding
  deterministically modeled system.  We also generalize our main
  result to some non-mass-action kinetics.
\end{abstract}


\setcounter{equation}{0}

\section{Introduction}
\label{sec:intro}

There are two commonly used models for chemical reaction systems:
discrete stochastic models in which the state of the system is a
vector giving the number of each molecular species, and continuous
deterministic models in which the state of the system is a vector
giving the concentration of each molecular species.  Discrete
stochastic models are typically used when the number of molecules of
each chemical species is low and the randomness inherent in the making
and breaking of chemical bonds is important.  Conversely,
deterministic models are used when there are large numbers of
molecules for each species and the behavior of the concentration of
each species is well approximated by a coupled set of ordinary
differential equations.

Typically, the goal in the study of discrete stochastic systems is to
either understand the evolution of the distribution of the state of
the system or to find the long term stationary distribution of the
system, which is the stochastic analog of an equilibrium point.  The
Kolmogorov forward equation (chemical master equation in the chemistry
literature) describes the evolution of the distribution and so work
has been done in trying to analyze or solve the forward equation for
certain classes of systems (\cite{Othmer2005}).  However, it is
typically an extremely difficult task to solve or even numerically
compute the solution to the forward equation for all but the simplest
of systems.  Therefore, simulation methods have been developed that
will generate sample paths so as to approximate the distribution of
the state via Monte Carlo methods.  These simulation methods include
algorithms that generate statistically exact (\cite{Anderson2007a,
  Gill76, Gill77, Gibson2000}) and approximate (\cite{Anderson2007b,
  AndGangKurtz2010, Gill2001, Cao2006}) sample paths.  On the other
hand, the continuous deterministic models, and in particular
mass-action systems with complex balancing states, have been analyzed
extensively in the mathematical chemistry literature, starting with
the works of Horn, Jackson, and Feinberg (\cite{Horn72, Horn73,
  HornJack72, Feinberg72}), and continuing with Feinberg's deficiency
theory in (\cite{FeinbergLec79, Feinberg87, Feinberg89,
  Feinberg95a}). Such models have a wide range of applications in the
physical sciences, and now they are beginning to play an important
role in systems biology (\cite{Craciun2006_PNAS, Gun2003,
  Sontag2001}).  Recent mathematical analysis of continuous
deterministic models has focused on their potential to admit multiple
equilibria (\cite{Craciun2005, Craciun2006}) and on dynamical
properties such as persistence and global stability (\cite{Sontag2001,
  Sontag2007, AndGlobal, AndCraciun2008, AndShiu2009}).

One of the major theorems pertaining to deterministic models of
chemical systems is the deficiency zero theorem of Feinberg
(\cite{Feinberg87, FeinbergLec79}).  The deficiency zero theorem
states that if the network of a system satisfies certain easily
checked properties, then within each compatibility class (invariant
manifold in which a solution is confined) there is precisely one
equilibrium with strictly positive components, and that equilibrium is
locally asymptotically stable (\cite{Feinberg87, FeinbergLec79}).  The
surprising aspect of the deficiency zero theorem is that the
assumptions of the theorem are completely related to the network of
the system whereas the conclusions of the theorem are related to the
dynamical properties of the system.  We will show in this paper that
if the conditions of the deficiency zero theorem hold on the network
of a stochastically modeled chemical system with quite general
kinetics, then there exists a product-form stationary distribution for
each closed, irreducible subset of the state space.  In fact, we will
show a stronger result: that a product-form stationary distribution
exists so long as there exists a complex balanced equilibrium for the
associated deterministically modeled system.  However, the equilibrium
values guaranteed to exist by the deficiency zero theorem are complex
balanced and so the conditions of that theorem are sufficient to
guarantee the existence of the product-form distribution.  Finally,
the main parameter of the stationary distribution will be shown to be
a complex balanced equilibrium value of the deterministically modeled
system.

Product-form stationary distributions play a central role in the
theory of queueing networks where the product-form property holds for
a large, naturally occurring class of models called Jackson networks
(see, for example, \cite{Kelly1979}, Chapter 3, and \cite{CY01},
Chapter 2) and a much larger class of quasi-reversible networks
(\cite{Kelly1979}, Chapter 3, \cite{CY01}, Chapter 4, \cite{Serf99},
Chapter 8).  Kelly, \cite{Kelly1979}, Section 8.5, recognizes the
possible existence of product-form stationary distributions for a
subclass of chemical reaction models and gives a condition for that
existence.  That condition is essentially the complex balance
condition described below, and our main result asserts that for any
mass-action chemical reaction model the conditions of the deficiency
zero theorem ensure that this condition holds.

The outline of the paper is as follows. In Section \ref{sec:networks}
we formally introduce chemical reaction networks.  In Section
\ref{sec:dyn} we develop both the stochastic and deterministic models
of chemical reaction systems.  Also in Section \ref{sec:dyn} we state
the deficiency zero theorem for deterministic systems and present two
theorems that are used in its proof and that will be of use to us.  In
Section \ref{sec:main_result_MA} we present the first of our main
results: that every closed, irreducible subset of the state space of a
stochastically modeled system with mass-action kinetics has a
product-form stationary distribution if the chemical network is weakly
reversible and has a deficiency of zero.  In Section
\ref{sec:examples} we present some examples of the use of this result.
In Section \ref{sec:result_Gen} we extend our main result to systems
with more general kinetics.

\setcounter{equation}{0}
\section{Chemical reaction networks}
\label{sec:networks}

Consider a system with $m$ chemical species, $\{S_1,\dots,S_m\}$,
undergoing a finite series of chemical reactions.  For the $k$th
reaction, denote by $\nu_k, \nu_k' \in \Z^m_{\ge 0}$ the vectors
representing the number of molecules of each species consumed and
created in one instance of that reaction, respectively.  We note that
if $\nu_k = \vec 0$ then the $k$th reaction represents an input to the
system, and if $\nu_k' = \vec 0$ then it represents an output.  Using
a slight abuse of notation, we associate each such $\nu_k$ (and
$\nu_k'$) with a linear combination of the species in which the
coefficient of $S_i$ is $\nu_{ik}$, the $i$th element of $\nu_k$.  For
example, if $\nu_k = [1, \ 2, \ 3]^T$ for a system consisting of three
species, we associate with $\nu_k$ the linear combination $S_1 + 2S_2
+ 3S_3$.  For $\nu_k = \vec 0$, we simply associate $\nu_k$ with
$\emptyset$.  Under this association, each $\nu_k$ (and $\nu_k'$) is
termed a {\em complex} of the system.  We denote any reaction by the
notation $\nu_k \to \nu_k'$, where $\nu_k$ is the source, or reactant,
complex and $\nu_k'$ is the product complex.  We note that each
complex may appear as both a source complex and a product complex in
the system.  The set of all complexes will be denoted by $\{\nu_k\} :=
\cup_k (\{\nu_k\} \cup \{\nu_k'\})$.

\begin{definition}
  Let $\S = \{S_i\}$, $\C = \{\nu_k\},$ and $\Re = \{\nu_k \to
  \nu_k'\}$ denote the sets of species, complexes, and reactions,
  respectively.  The triple $\{\S, \C, \Re\}$ is called a {\em
    chemical reaction network}.
  \label{def:crn}
\end{definition}

The structure of chemical reaction networks plays a central role in
both the study of stochastically and deterministically modeled
systems.  As alluded to in the Introduction, it will be conditions on
the network of a system that guarantee certain dynamical properties
for both models.  Therefore, the remainder of this section consists of
definitions related to chemical networks that will be used throughout
the paper.

\begin{definition}
  A chemical reaction network, $\{\S,\C,\Re \}$, is called {\em weakly
    reversible} if for any reaction $\nu_k \to \nu_k'$, there is a
  sequence of directed reactions beginning with $\nu_k'$ as a source
  complex and ending with $\nu_k$ as a product complex.  That is,
  there exist complexes $\nu_1,\dots,\nu_r$ such that $\nu_k' \to
  \nu_1, \nu_1 \to \nu_2, \dots, \nu_r \to \nu_k \in \Re$.  A network
  is called {\em reversible} if $\nu_k' \to \nu_k \in \Re$ whenever
  $\nu_k \to \nu_k' \in \Re$.
  \label{def:WR}
\end{definition}

\begin{remark} The definition of a reversible network given in
  Definition \ref{def:WR} is distinct from the notion of a reversible
  stochastic process.  However, in Section \ref{sec:rev} we point out
  a connection between the two concepts for systems that are detailed
  balanced.
\end{remark}

To each reaction network, $\{\S,\C,\Re\}$, there is a unique, directed
graph constructed in the following manner.  The nodes of the graph are
the complexes, $\C$.  A directed edge is then placed from complex
$\nu_k$ to complex $\nu_k'$ if and only if $\nu_k \to \nu_k' \in \Re$.
Each connected component of the resulting graph is termed a {\em
  linkage class} of the graph.  We denote the number of linkage
classes by $\ell$.  It is easy to see that a chemical reaction network
is weakly reversible if and only if each of the linkage classes of its
graph is strongly connected.

\begin{definition}
  $S = \hbox{span}_{\{\nu_k \to \nu_k' \in \Re\}}\{\nu_k' - \nu_k\}$
  is the {\em stoichiometric subspace} of the network. For $c \in
  \R^m$ we say $c + S$ and $(c + S) \cap \R^m_{>0}$ are the {\em
    stoichiometric compatibility classes} and {\em positive
    stoichiometric compatibility classes} of the network,
  respectively.  Denote $\hbox{dim}(S) = s$.
\end{definition}

It is simple to show that for both stochastic and deterministic
models, the state of the system remains within a single stoichiometric
compatibility class for all time, assuming that one starts in that
class.  This fact is important because it changes the types of
questions that are reasonable to ask about a given system.  For
example, unless there is only one stoichiometric compatibility class,
and so $S = \R^m$, the correct question is not whether there is a
unique fixed point for a given deterministic system.  Instead, the
correct question is \textit{whether within each stoichiometric
  compatibility class} there is a unique fixed point.  Analogously,
for stochastically modeled systems it is typically of interest to
compute stationary distributions for each closed, irreducible subset
of the state space (each contained within a stoichiometric
compatibility class) with the precise subset being determined by
initial conditions.

The final definition of this section is that of the {\em deficiency}
of a network (\cite{FeinbergLec79}).  It is not a difficult exercise
to show that the deficiency of a network is always greater than or
equal to zero.
\begin{definition}
  The {\em deficiency} of a chemical reaction network,
  $\{\S,\C,\Re\}$, is $\delta = |\C| - \ell - s$, where $|\C|$ is the
  number of complexes, $\ell$ is the number of linkage classes of the
  network graph, and $s$ is the dimension of the stoichiometric
  subspace of the network.
\end{definition}
While the deficiency is, by definition, only a property of the
network, we will see in Sections \ref{sec:def0},
\ref{sec:main_result_MA}, and \ref{sec:result_Gen} that a deficiency
of zero has implications for the long-time dynamics of both
deterministic and stochastic models of chemical reaction systems.

\setcounter{equation}{0}
\section{Dynamical models}
\label{sec:dyn}

The notion of a chemical reaction network is the same for both
stochastic and deterministic systems and the choice of whether to
model the evolution of the state of the system stochastically or
deterministically is made based upon the details of the specific
chemical or biological problem at hand.  Typically if the number of
molecules is low, a stochastic model is used, and if the number of
molecules is high, a deterministic model is used.  For cases between
the two extremes a diffusion approximation can be used or, for cases
in which the system contains multiple scales, pieces of the reaction
network can be modeled stochastically, while others can be modeled
deterministically (or, more accurately, absolutely continuously with
respect to time).  See, for example, \cite{Ball06} and Section
\ref{sec:multiscale}.

\subsection{Stochastic models}

The simplest stochastic model for a chemical network $\{\S, \C, \Re\}$
treats the system as a continuous time Markov chain whose state $X \in
\Z^m_{\ge 0}$ is a vector giving the number of molecules of each
species present with each reaction modeled as a possible transition
for the state.  We assume a finite number of reactions.  The model for
the $k$th reaction, $\nu_k \to \nu_k'$, is determined by the vector of
inputs, $\nu_k$, specifying the number of molecules of each chemical
species that are consumed in the reaction, the vector of outputs,
$\nu_k'$, specifying the number of molecules of each species that are
created in the reaction, and a function of the state, $\lambda_k(X)$,
that gives the rate at which the reaction occurs.  Specifically, if
the $k$th reaction occurs at time $t$, the new state becomes
\begin{equation*}
  X(t) = X(t-) + \nu_k' - \nu_k.
\end{equation*}
Let $R_k(t)$ denote the number of times that the $k$th reaction occurs
by time $t$.  Then the state of the system at time $t$ can be written
as
\begin{equation}
  X(t) = X(0) + \sum_{k} R_k(t)(\nu_k' - \nu_k),
  \label{eq:stoch_model}
\end{equation}
where we have summed over the reactions.  The process $R_k$ is a
counting process with intensity $\lambda_k(X(t))$ (called the {\em
  propensity} in the chemistry literature) and can be written as
\begin{equation}
  R_k(t) = Y_k\left(\int_0^t \lambda_k(X(s))ds \right),
  \label{eq:Rk}
\end{equation}
where the $Y_k$ are independent, unit-rate Poisson processes
(\cite{Kurtz78}, \cite{Kurtz86} Ch. 11). Note that
\eqref{eq:stoch_model} and \eqref{eq:Rk} give a system of stochastic
equations that uniquely determines $X$ up to $\sup \{t \, :\, \sum_k
R_k(t) < \infty\}$.  The generator for the Markov chain is the
operator, $A$, defined by
\begin{equation}
  Af(x) = \sum_{k} \lambda_k(x)(f(x + \nu_k' - \nu_k) - f(x)),
  \label{eq:generator}
\end{equation}
where $f$ is any function defined on the state space.

A commonly chosen form for the intensity functions $\lambda_k$ is that
of stochastic mass-action, which says that for $x \in \Z^m_{\ge 0}$
the rate of the $k$th reaction should be given by
\begin{equation}
  \lambda_k(x)  = \kappa_k \left(\prod_{\ell =1}^m \nu_{\ell k}!
  \right) \binom{x}{\nu_k} = \kappa_k \prod_{\ell = 1}^m 
  \frac{x_{\ell}!}{(x_{\ell} - \nu_{\ell k})!} 1_{\{x_{\ell} \ge \nu_{\ell k}\}},
  \label{eq:stoch_MA}
\end{equation}
for some constant $\kappa_k$, where we adopt the convention that $0! =
1$.  Note that the rate \eqref{eq:stoch_MA} is proportional to the
number of distinct subsets of the molecules present that can form the
inputs for the reaction.  Intuitively, this assumption reflects the
idea that the system is {\em well-stirred} in the sense that all
molecules are equally likely to be at any location at any time.  For
concreteness, we will assume that the intensity functions satisfy
\eqref{eq:stoch_MA} throughout most of the paper. In Section
\ref{sec:result_Gen} we will generalize our results to systems with
more general kinetics.

A probability distribution $\{\pi(x)\}$ is a stationary distribution
for the chain if
\begin{equation*}
  \sum_{x} \pi(x)Af(x) = 0
\end{equation*}
for a sufficiently large class of functions $f$ or, taking $f(y) =
\mathbf{1}_{x}(y)$ and using equation \eqref{eq:generator}, if
\begin{equation}
  \sum_{k} \pi(x - \nu_k' + \nu_k)\lambda_k(x - \nu_k' + \nu_k) =
  \pi(x)\sum_k \lambda_k(x)
  \label{eq:stationary2}
\end{equation}
for all $x$ in the state space.  If the network is weakly reversible,
then the state space of the Markov chain is a union of closed,
irreducible communicating classes.  (This fact follows because if the
Markov chain can proceed from state $x$ to state $y$ via a sequence of
reactions, weak reversibility of the network implies those reactions
can be ``undone'' in reverse sequential order by another sequence of
reactions.)  Also, each closed, irreducible communicating class is
either finite or countable.  Therefore, if a stationary distribution
with support on a single communicating class exists it is unique and
\begin{equation*}
  \lim_{t \to \infty} P(X(t) = x \ | \ X(0) = y) = \pi(x),
\end{equation*}
for all $x,y$ in that communicating class.  Thus, the stationary
distribution gives the long-term behavior of the system.

Solving equation \eqref{eq:stationary2} is in general a formidable
task.  However, in Section \ref{sec:main_result_MA} we will do so if
the network is weakly reversible, has a deficiency of zero, and if the
rate functions $\lambda_k(x)$ satisfy mass-action kinetics,
\eqref{eq:stoch_MA}.  We will also show that the stationary
distribution is of product form.  More specifically, we will show that
for each communicating class there exists a $c \in \mathbb{R}^m_{>0}$
and a normalizing constant $M>0$ such that
\begin{equation*}
  \pi(x) = M \prod_{i=1}^m \pi_i(x_i) :=  M \prod_{i=1}^m
  \frac{c_i^{x_i}}{x_i!}
\end{equation*}
satisfies equation \eqref{eq:stationary2}.  The $c_i$ in the
definition of $\pi_i$ will be shown to be the $i$th component of an
equilibrium value of the analogous deterministic system described in
the next section.  In Section \ref{sec:result_Gen} we will solve
\eqref{eq:stationary2} for more general kinetics.

\subsection{Deterministic models and the deficiency zero theorem}
\label{sec:def0}

Under an appropriate scaling limit (see Section \ref{sec:scaling})
the continuous time Markov chain \eqref{eq:stoch_model},
\eqref{eq:Rk}, \eqref{eq:stoch_MA} becomes
\begin{equation}
  x(t) = x(0) + \sum_{k} \left(\int_0^t f_k(x(s)) ds \right) (\nu_k' -
  \nu_k) := x(0) + \int_0^t f(x(s)) ds,
  \label{eq:cont}
\end{equation}
where the last equality is a definition and
\begin{equation}
  f_{k}(x) = \ \kappa_{k}x_1^{\nu_{1k}} x_2^{\nu_{2k}} \cdots x_m^{\nu_{mk}},
  \label{eq:massaction}
\end{equation}
where we use the convention $0^0 = 1$.  We say that the deterministic
system \eqref{eq:cont} has {\em mass-action kinetics} if the rate
functions $f_k$ have the form \eqref{eq:massaction}.  The proof of the
following theorem by Feinberg can be found in \cite{FeinbergLec79} or
\cite{Feinberg95a}.  We note that the full statement of the deficiency
zero theorem actually says more than what is given below and the
interested reader is encouraged to see the original work.

\begin{theorem}[The Deficiency Zero Theorem]
  Consider a weakly reversible, deficiency zero chemical reaction
  network $\{\S, \C, \Re\}$ with dynamics given by
  \eqref{eq:cont}-\eqref{eq:massaction}.  Then for any choice of rate
  constants $\{\kappa_k\}$, within each positive stoichiometric
  compatibility class there is precisely one equilibrium value, and
  that equilibrium value is locally asymptotically stable relative to
  its compatibility class.
  \label{thm:def0}
\end{theorem}

The dynamics of the system \eqref{eq:cont}-\eqref{eq:massaction} take
place in $\R^m_{\ge 0}$.  However, to prove the deficiency zero
theorem it turns out to be more appropriate to work in {\em complex
  space}, denoted $\R^{\C}$, which we will describe now.  For any $U
\subseteq \C$ let $\omega_U: \C \to \{0,1\}$ denote the indicator
function $\omega_U(\nu_k) = \mathbf{1}_{\{\nu_k \in U\}}$. Complex
space is defined to be the vector space with basis $\{\omega_{\nu_k} \
| \ \nu_k \in \C\}$, where we have denoted $\omega_{\{\nu_k\}}$ by
$\omega_{\nu_k}$.

If $u$ is a vector with nonnegative integer components and $w$ is a
vector with nonnegative real components, then let $u! = \prod_i u_i!$
and $w^u = \prod_i w_i^{u_i}$, where we interpret $0^0 = 1$ and $0! =
1$.  Let $\Psi: \R^{m} \to \R^{\C}$ and $A_{\kappa}: \R^{\C} \to
\R^{\C}$ be defined by:
\begin{align*}
  \Psi(x) &= \sum_{\nu_k \in \C} x^{\nu_k}\omega_{\nu_k}\\
  A_{\kappa}(y) &= \sum_{\nu_k \to \nu_k' \in \Re} \kappa_k
  y_{\nu_k}(\omega_{\nu_k'} - \omega_{\nu_k}),
\end{align*}
where the subscript $\kappa$ of $A_{\kappa}$ denotes the choice of
rate constants for the system.  Let $Y: \R^{\C} \to \R^{m}$ be the
linear map whose action on the basis elements $\{\omega_{\nu_k}\}$ is
defined by $Y(\omega_{\nu_k}) = \nu_k$.  Then equations
\eqref{eq:cont}-\eqref{eq:massaction} can be written as the coupled
set of ordinary differential equations
\begin{equation*}
  \dot x(t) = f(x(t)) = Y(A_{\kappa}(\Psi(x(t)))).
\end{equation*}
Therefore, in order to show that a value $c$ is an equilibrium of the
system, it is sufficient to show that $A_{\kappa}(\Psi(c)) = 0$, which
is an explicit system of equations for $c$.  In particular,
$A_k(\Psi(c)) = 0$ if and only if for each $z \in \C$
\begin{equation}
  \sum_{\{k : \nu_k' = z \}} \kappa_k c^{\nu_k} = \sum_{ \{ k:\nu_k=z
    \} } \kappa_k c^{\nu_k}, 
  \label{eq:complex_balanced}
\end{equation}
where the sum on the left is over reactions for which $z$ is the
product complex and the sum on the right is over reactions for which
$z$ is the source complex.

The following has been shown in \cite{HornJack72} and
\cite{FeinbergLec79} (see also \cite{Gun2003}).
\begin{theorem}
  Let $\{\S, \C, \Re\}$ be a chemical reaction network with dynamics
  given by \eqref{eq:cont}-\eqref{eq:massaction} for some choice of
  rate constants, $\{\kappa_k\}$.  Suppose there exists a $c \in
  \R^m_{>0}$ for which $A_{\kappa}(\Psi(c)) = 0$, then the following
  hold:
  \begin{enumerate}
  \item The network is weakly reversible.
  \item Every equilibrium point with strictly positive components, $x
    \in \R^m_{>0}$ with $f(x) = 0$, satisfies $A_{\kappa}(\Psi(x)) =
    0$.
  \item If $Z = \{x \in \R^m_{>0} \ | \ f(x) = 0\}$, then $\ln Z :=
    \{y \in \mathbb{R}^m \ | \ \exists \ x \in Z \hbox{ and } y_i =
    \ln(x_i)\}$ is a coset of $S^{\perp}$, the perpendicular
    complement of $S$.  That is, there is a $k \in \mathbb{R}^m$ such
    that $\ln Z = \{w \in \mathbb{R}^m \ | \ w = k + u \hbox{ for some
    } u \in S^{\perp}\}$.
  \item There is one, and only one, equilibrium point in each positive
    stoichiometric compatibility class.
  \item Each equilibrium point of a positive stoichiometric
    compatibility class is locally asymptotically stable relative to
    its stoichiometric compatibility class.
  \end{enumerate}
  \label{thm:def0_main}
\end{theorem}
Thus, after a choice of rate constants has been made, the conclusions
of the deficiency zero theorem pertaining to the existence and
asymptotic stability of equilibria (points $4.$ and $5.$ of Theorem
\ref{thm:def0_main}) hold so long as there exists at least one $c \in
\R^m_{>0}$ such that $A_{\kappa}(\Psi(c)) = 0$.  The condition that
the system has a deficiency of zero only plays a role in showing that
there does exist such a $c \in \R^m_{>0}$.  A proof of the following
can be found in \cite{FeinbergLec79}, \cite{Feinberg87}, or
\cite{Feinberg95a}.

\begin{theorem}
  Let $\{\S, \C, \Re\}$ be a chemical reaction network with dynamics
  given by \eqref{eq:cont}-\eqref{eq:massaction} for some choice of
  rate constants, $\{\kappa_k\}$.  If the network has a deficiency of
  zero, then there exists a $c \in \R^m_{>0}$ such that
  $A_{\kappa}(\Psi(c)) = 0$ if and only if the network is weakly
  reversible.
  \label{thm:def0_def}
\end{theorem}

A chemical reaction network with deterministic mass-action kinetics
(and a choice of rate constants) that admits a $c$ for which
$A_{\kappa}(\Psi(c)) = 0$ is called \textit{complex balanced} in the
literature. The second conclusion of Theorem \ref{thm:def0_main}
demonstrates why this notation is appropriate.  The equivalent
representation given by \eqref{eq:complex_balanced} shows the origin
of this terminology.  The surprising aspect of the deficiency zero
theorem is that it gives simple and checkable sufficient conditions on
the network structure alone that guarantee that a system is complex
balanced for any choice of rate constants.  We will see in the
following sections that the main results of this paper have the same
property: product-form stationary distributions exist for all
stochastic systems that are complex balanced when viewed as
deterministic systems, and $\delta = 0$ is a sufficient condition to
guarantee this for weakly reversible networks.

\setcounter{equation}{0}
\section{Main result for mass-action systems}
\label{sec:main_result_MA}

The collection of stationary distributions for a countable state space
Markov chain is convex. The extremal distributions correspond to the
closed, irreducible subsets of the state space; that is, every
stationary distribution can be written as
\begin{equation}
  \pi = \sum_{\Gamma} \alpha_{\Gamma} \pi_{\Gamma},
  \label{eq:pi_gamma}
\end{equation}
where $\alpha_{\Gamma} \ge 0$, $\sum_{\Gamma} \alpha_{\Gamma} = 1$,
and the sums are over the closed, irreducible subsets $\Gamma$ of the
state space.  Here $\pi_{\Gamma}$ is the unique stationary
distribution satisfying $\pi_{\Gamma}(\Gamma) = 1$.

We now state and prove our main result for systems with mass-action
kinetics.

\begin{theorem}
  Let $\{\S, \C, \Re\}$ be a chemical reaction network and let
  $\{\kappa_k\}$ be a choice of rate constants.  Suppose that, modeled
  deterministically, the system is complex balanced with complex
  balanced equilibrium $c \in \R^m_{>0}$.  Then the stochastically
  modeled system with intensities \eqref{eq:stoch_MA} has a stationary
  distribution consisting of the product of Poisson distributions,
  \begin{equation}
    \pi(x) = \prod_{i=1}^m \frac{c_i^{x_i}}{x_i!}e^{-c_i}, \qquad x \in
    \Z^m_{\ge 0}.
    \label{eq:stationary_product}
  \end{equation}
  If $\Z^m_{\ge 0}$ is irreducible, then \eqref{eq:stationary_product}
  is the unique stationary distribution, whereas if $\Z^m_{\ge 0}$ is
  not irreducible then the $\pi_{\Gamma}$ of equation
  \eqref{eq:pi_gamma} are given by the product-form stationary
  distributions
  \begin{equation*}
    \pi_{\Gamma}(x) = M_{\Gamma} \prod_{i = 1}^m
    \frac{c_i^{x_i}}{x_i!}, \qquad x \in \Gamma,
  \end{equation*}
  and $\pi_{\Gamma}(x) = 0$ otherwise, where $M_{\Gamma}$ is a
  positive normalizing constant.
  \label{thm:prodform_main}
\end{theorem}

\begin{proof}
  
  Let $\pi$ satisfy \eqref{eq:stationary_product} where $c\in
  \R^m_{>0}$ satisfies $A_{\kappa}(\Psi(c)) = 0$.  We will show that
  $\pi$ is stationary by verifying that equation
  \eqref{eq:stationary2} holds for all $x \in \Z^{m}_{\ge 0}$.
  Plugging $\pi$ and \eqref{eq:stoch_MA} into equation
  \eqref{eq:stationary2} and simplifying yields
  \begin{equation}
    \sum_{k} \kappa_k c^{ \nu_k - \nu_k'} \frac{1}{(x -
      \nu_k')!} \prod_{\ell = 1}^m 1_{\{x_{\ell} \ge \nu_{\ell k}' \}}
    = \sum_k \kappa_k \frac{1}{(x - \nu_k)!} \prod_{\ell = 1}^m
    1_{\{x_{\ell} \ge \nu_{\ell k} \}}. 
    \label{eq:stat3}
  \end{equation}
  Equation \eqref{eq:stat3} will be satisfied if for each complex $z
  \in \C$,
  \begin{equation}
    \sum_{\{k:\nu_k'=z\}} \kappa_k c^{\nu_k-z} \frac{1}{(x - z)!}
    \prod_{\ell = 1}^m 1_{\{x_{\ell} \ge z_{\ell} \}}=
    \sum_{\{k:\nu_k=z\}} \kappa_k \frac{1}{(x - z)!} \prod_{\ell =
      1}^m 1_{\{x_{\ell} \ge z_{\ell} \}}, 
    \label{eq:stat4}
  \end{equation}
  where the sum on the left is over reactions for which $z$ is the
  product complex and the sum on the right is over reactions for which
  $z$ is the source complex.  The complex $z$ is fixed in the above
  equation, and so \eqref{eq:stat4} is equivalent to
  \eqref{eq:complex_balanced}, which is equivalent to
  $A_{\kappa}(\Psi(c)) = 0$.

  To complete the proof, one need only observe that the normalized
  restriction of $\pi$ to any closed, irreducible subset $\Gamma$ must
  also be a stationary distribution.
\end{proof}

The following theorem gives simple and checkable conditions that
guarantee the existence of a product-form stationary distribution of
the form \eqref{eq:stationary_product}.

\begin{theorem}
  Let $\{\S, \C, \Re\}$ be a chemical reaction network that has a
  deficiency of zero and is weakly reversible.  Then for any choice of
  rate constants $\{\kappa_k\}$ the stochastically modeled system with
  intensities \eqref{eq:stoch_MA} has a stationary distribution
  consisting of the product of Poisson distributions,
  \begin{equation*}
    \pi(x) = \prod_{i=1}^m \frac{c_i^{x_i}}{x_i!} e^{-c_i}, \qquad x \in
    \Z^m_{\ge 0}, 
  \end{equation*}
  where $c$ is an equilibrium value for the deterministic system
  \eqref{eq:cont}-\eqref{eq:massaction}, which is guaranteed to exist
  and be complex balanced by Theorems
  \ref{thm:def0}-\ref{thm:def0_def}.  If $\Z^m_{\ge 0}$ is
  irreducible, then $\pi$ is the unique stationary distribution,
  whereas if $\Z^m_{\ge 0}$ is not irreducible then the $\pi_{\Gamma}$
  of equation \eqref{eq:pi_gamma} are given by the product-form
  stationary distributions
  \begin{equation*}
    \pi_{\Gamma}(x) = M_{\Gamma} \prod_{i = 1}^m
    \frac{c_i^{x_i}}{x_i!}, \qquad x \in \Gamma, 
  \end{equation*}
  and $\pi_{\Gamma}(x) = 0$ otherwise, where $M_{\Gamma}$ is a
  positive normalizing constant.
  \label{thm:prodform_def0}
\end{theorem}

\begin{proof}
  This is a direct result of Theorems \ref{thm:def0_def} and
  \ref{thm:prodform_main}.
\end{proof}

We remark that Theorems \ref{thm:prodform_main} and
\ref{thm:prodform_def0} give sufficient conditions under which
$\Z^m_{\ge 0}$ being irreducible guarantees that when in
distributional equilibrium the species numbers: (a) are independent
and (b) have Poisson distributions.  We return to this point in
Examples \ref{example:first_order} and \ref{ex:first_enzyme}.

\subsection{The classical scaling}
\label{sec:scaling}

Defining $|\nu_k| =\sum_i \nu_{ik}$ and letting $V$ be a scaling
parameter usually taken to be the volume of the system times
Avogadro's number, it is reasonable to scale the rate constants of the stochastic model with the volume like
\begin{equation}
 \kappa_k = \frac{\hat \kappa_k}{V^{|\nu_k| - 1}},
 \label{eq:kappas}
\end{equation}
for some $\hat \kappa_k >0$.
This follows by considering the probability of a particular set of $|\nu_k|$ molecules finding each other in a volume proportional to $V$ in a time interval $[t,t+\Delta t)$.  In this case, the intensity functions become
\begin{equation}
  \lambda_k^V(x) =  \frac{\hat{\kappa}_k}{V^{|\nu_k| - 1}}
  (\prod_{i} \nu_{ik}!) \binom{x}{\nu_{k}} = V\hat{\kappa}_k
  \frac{1}{V^{|\nu_k|}}\prod_{i} \frac{x_{i}!}{(x_{i} - \nu_{ik
    })!}.    
    \label{eq:classic_scaling}
\end{equation}
   Since $V$ is the volume times Avogadro's number and $x$ gives the
number of molecules of each species present, $c = V^{-1}x$ gives the
concentrations in moles per unit volume. With this scaling and a large
volume limit
\begin{equation}
  \lambda_k^V(x) \approx V \hat{\kappa}_k \prod_i c_i^{\nu_{ik}} = V \hat \kappa_k c^{\nu_k} \equiv
  V \hat{\lambda}_k(c). 
  \label{eq:first_approx}
\end{equation}
Since the law of large numbers for the Poisson process implies
$V^{-1}Y_k(Nu)\approx u$, \eqref{eq:Rk} and \eqref{eq:first_approx}, together with the assumption that $X(0) =VC(0)$ for some $C(0)\in \R^m_{>0}$, imply
\begin{equation*}
  C(t) = V^{-1} X(t) \approx C(0) + \sum_k \int_0^t \hat{\kappa}_k
  C(s)^{\nu_{k}} ds \, (\nu_{k}' - \nu_k),
\end{equation*}
which in the large volume limit gives the classical deterministic law
of mass action detailed in Section \ref{sec:def0}.   For a precise formulation of the above scaling
argument, termed the ``classical scaling,'' see \cite{Kurtz72,
  Kurtz78, KurtzPop81}. 
  
Because the above scaling is the natural relationship between the stochastic and deterministic models of chemical reaction networks, we expect to be able to generalize Theorem \ref{thm:prodform_main} to this setting.

\begin{theorem}
 Let $\{\S, \C, \Re\}$ be a chemical reaction network.   Suppose that, modeled
  deterministically with rate constants $\{\hat \kappa_k\}$, the system is complex balanced with complex
  balanced equilibrium $c \in \R^m_{>0}$.  For some $V>0$, let $\{\kappa_k\}$ be related to $\{\hat \kappa_k\}$ via \eqref{eq:kappas}.  Then the stochastically
  modeled system with intensities \eqref{eq:stoch_MA} and rate constants $\{\kappa_k\}$ has a stationary distribution consisting of the product of Poisson distributions,
  \begin{equation*}
    \pi(x) = \prod_{i=1}^m \frac{(Vc_i)^{x_i}}{x_i!}e^{-Vc_i}, \qquad x \in
    \Z^m_{\ge 0}.
  \end{equation*}
  If $\Z^m_{\ge 0}$ is irreducible, then \eqref{eq:stationary_product}
  is the unique stationary distribution, whereas if $\Z^m_{\ge 0}$ is
  not irreducible then the $\pi_{\Gamma}$ of equation
  \eqref{eq:pi_gamma} are given by the product-form stationary
  distributions
  \begin{equation*}
    \pi_{\Gamma}(x) = M_{\Gamma} \prod_{i = 1}^m
    \frac{(Vc_i)^{x_i}}{x_i!}, \qquad x \in \Gamma,
  \end{equation*}
  and $\pi_{\Gamma}(x) = 0$ otherwise, where $M_{\Gamma}$ is a
  positive normalizing constant.
  \label{thm:rate_constants}
\end{theorem}

\begin{proof}
The proof is similar to before, and now consists of making sure the $V$'s cancel in an appropriate manner.  The details are omitted. 
\end{proof}
We see that Theorem \ref{thm:prodform_main} follows from Theorem \ref{thm:rate_constants} by taking $V = 1$.  Theorem \ref{thm:prodform_def0} generalizes in the obvious way.

\subsection{Reversibility and detail balance}
\label{sec:rev}

An equilibrium value, $c \in \R^m_{>0}$, for a reversible, in the
sense of Definition \ref{def:WR}, chemical reaction network with
deterministic mass-action kinetics is called \textit{detailed
  balanced} if for each pair of reversible reactions, $\nu_k
\rightleftarrows \nu_k'$, we have
\begin{equation}
  \kappa_k c^{\nu_k} = \kappa_k' c^{\nu_k'},
  \label{eq:det_bal1}
\end{equation}
where $\kappa_k,\kappa_k'$ are the rate constants for the reactions
$\nu_k \to \nu_k', \nu_k' \to \nu_k$, respectively.  In
\cite{Feinberg89}, page 1820, Feinberg shows that if one positive
equilibrium is detailed balanced then they all are; a result similar
to the second conclusion of Theorem \ref{thm:def0_main} for complex
balanced systems. A reversible chemical reaction system with
deterministic mass action kinetics is therefore called
\textit{detailed balanced} if it admits one detailed balanced
equilibrium. It is immediate that any system that is detailed balanced
is also complex balanced.  The fact that a product-form stationary
distribution of the form \eqref{eq:stationary_product} exists for the
stochastic systems whose deterministic analogs are detailed balanced
is well-known.  See, for example, \cite{Whittle86}.  Theorems
\ref{thm:prodform_main} and \ref{thm:prodform_def0} can therefore be
viewed as an extension of that result.  However, more can be said in the case when the deterministic system is
detailed balanced, and which we include here for completeness (no originality is being claimed).

  As mentioned in the remark following Definition
\ref{def:WR}, the term ``reversible'' has a meaning in the context of
stochastic processes that differs from that of Definition
\ref{def:WR}.   Before defining this, we need the concept of a transition rate.  For any continuous time Markov chain with state space
$\Gamma$, the \textit{transition rate} from $x\in \Gamma$ to $y \in \Gamma$ (with $x\ne y$) is a non-negative number $\alpha(x,y)$ satisfying 
\begin{equation*}
  P(X(t + \Delta t) = y\ | \ X(t) = x) = \alpha(x,y)\Delta t + o(\Delta t).
\end{equation*}
Thus, in the context of this paper, if $y = x + \nu_k' - \nu_k$ for some $k$, then
$\alpha(x,y) = \lambda_k(x)$, and zero otherwise.


\begin{definition}
  A continuous time Markov chain $X(t)$ with transition rates
  $\alpha(x,y)$ is {\em reversible with respect to the distribution
    $\pi$} if for all $x,y$ in the state space $\Gamma$
  \begin{equation}
      \pi(x) \alpha(x,y) = \pi(y) \alpha(y,x).
      \label{eq:reversible}
  \end{equation}
  \label{def:rev}
\end{definition}

It is simple to see (by summing both sides of \eqref{eq:reversible}
with respect to $y$ over $\Gamma$), that $\pi$ must be a
stationary distribution for the process.  A stationary distribution
satisfying \eqref{eq:reversible} is even called \textit{detailed
  balanced} in the probability literature.  The following is proved in \cite{Whittle86}, Chapter 7.

\begin{theorem}
  Let $\{\S, \C, \Re\}$ be a reversible$^{2}$ chemical reaction
  network with rate constants $\{\kappa_k\}$.  Then the
  deterministically modeled system with mass-action kinetics has a
  detailed balanced equilibrium if and only if the stochastically
  modeled system with intensities \eqref{eq:stoch_MA} is reversible
  with respect to its stationary distribution.$^3$
  \label{thm:reversibility}
\end{theorem}

Succinctly, this theorem says that reversibility and detailed balanced in the deterministic setting is equivalent to reversible (and, hence, detailed balanced) in the stochastic setting.

\footnotetext[2]{In the sense of Definition \ref{def:WR}.}
\footnotetext[3]{in the sense of Definition \ref{def:rev}.}

\subsection{Non-uniqueness of $c$}

For stochastically modeled chemical reaction systems any irreducible
subset of the state space, $\Gamma$, is contained within $(y + S) \cap
\Z^m_{\ge 0}$ for some $y \in \R^m_{\ge 0}$.  Therefore, each $\Gamma$
is associated with a stoichiometric compatibility class.  For weakly
reversible systems with a deficiency of zero, Theorems
\ref{thm:def0_main} and \ref{thm:def0_def} guarantee that each such
stoichiometric compatibility class has an associated equilibrium value
for which $A_{\kappa}(\Psi(c)) = 0$.  However, neither Theorem
\ref{thm:prodform_main} nor Theorem \ref{thm:prodform_def0} makes the
requirement that the equilibrium value used in the product-form
stationary measure $\pi_{\Gamma}(\cdot)$ be contained within the
stoichiometric compatibility class associated with $\Gamma$.
Therefore we see that one such $c$ can be used to construct a
product-form stationary distribution for every closed, irreducible
subset.  Conversely, for a given irreducible subset $\Gamma$ any
positive equilibrium value of the system
\eqref{eq:cont}-\eqref{eq:massaction} can be used to construct
$\pi_{\Gamma}(\cdot)$.  This fact seems to be contrary to the
uniqueness of the stationary distribution; however, it can be
understood through the third conclusion of Theorem \ref{thm:def0_main}
as follows.

Let $\Gamma$ be a closed, irreducible subset of the state space with
associated positive stoichiometric compatibility class $(y + S) \cap
\Z^m_{\ge 0}$, and let $c_1, c_2 \in \R^m_{>0}$ be such that
$A_{\kappa}(\Psi(c_1)) = A_{\kappa}(\Psi(c_2)) = 0$.  For $i \in
\{1,2\}$ and $x \in \Gamma$, let $\pi_i(x) = M_i c_i^x/x!$, where
$M_1$ and $M_2$ are normalizing constants.  Then for each $x \in
\Gamma$
\begin{equation*}
  \frac{\pi_1(x)}{\pi_2(x)} = \frac{M_1 c_1^x}{x!}\frac{x!}{M_2c_2^x}
  = \frac{M_1}{M_2} \frac{c_1^x}{c_2^x}.
\end{equation*}
For any vector $u$, we define $(\ln(u))_i = \ln(u_i)$.  Then for $x
\in \Gamma \subset y + S$
\begin{equation}
  \frac{c_1^x}{c_2^x} = e^{x \cdot \left(\ln c_1 - \ln c_2\right)}
  = e^{y \cdot \left(\ln c_1 - \ln c_2\right)}
  = \frac{c_1^y}{c_2^y},
\label{eq:equilib_equiv1}
\end{equation}
where the second equality follows from the third conclusion of Theorem
\ref{thm:def0_main}.  Therefore,
\begin{equation}
  \frac{\pi_1(x)}{\pi_2(x)} = \frac{M_1}{M_2}\frac{c_1^y}{c_2^y}.
\label{eq:equilib_equiv2}
\end{equation}
Finally,
\begin{align*}
  1 &= \left(M_1\sum_{x \in \Gamma} c_1^x/x!\right) / \left(M_2\sum_{x
      \in \Gamma} c_2^x/x!\right)\\
  &= \frac{M_1}{M_2} \left(\frac{c_1^y}{c_2^y}\sum_{x \in \Gamma}
    c_2^x/x!\right) / \left(\sum_{x \in \Gamma} c_2^x/x!\right)\\
  &= \frac{\pi_1(x)}{\pi_2(x)},
\end{align*}
where the second equality follows from equation
\eqref{eq:equilib_equiv1} and the third equality follows from equation
\eqref{eq:equilib_equiv2}.  We therefore see that the stationary
measure is independent of the choice of $c$, as expected.

\setcounter{equation}{0}
\section{Examples}
\label{sec:examples}

Our first example points out that the existence of a product-form
stationary distribution for the closed, irreducible subsets of the
state space does not necessarily imply independence of the species
numbers.

\begin{example}\textbf{(Non-independence of species numbers)}
  Consider the simple reversible system
  \begin{equation*}
    S_1 \begin{array}{c}
      k_1 \\
      \rightleftarrows \\
      k_2
    \end{array} S_2,
  \end{equation*}
  where $k_1$ and $k_2$ are nonzero rate constants.  We suppose that
  $X_1(0) + X_2(0) = N$, and so $X_1(t) + X_2(t) = N$ for all $t$.
  This system has two complexes, one linkage class, and the dimension
  of the stoichiometric compatibility class is one.  Therefore it has
  a deficiency of zero.  Since it is also weakly reversible, our
  results hold.  An equilibrium to the system that satisfies the
  complex balance equation is
  \begin{equation*}
    c = \left(\frac{k_2}{k_1 + k_2}, \frac{k_1}{k_1 + k_2}\right),
  \end{equation*}
  and the product-form stationary distribution for the system is
  \begin{equation*}
    \pi(x) = M \frac{c_1^{x_1}}{x_1!} \frac{c_2^{x_2}}{x_2!},
  \end{equation*}
  where $M>0$ is a normalizing constant. Using that $X_1(t) + X_2(t) =
  N$ for all $t$ yields
  \begin{equation*}
    \pi_1(x_1) = M  \frac{c_1^{x_1}}{x_1!} \frac{c_2^{N - x_1}}{(N -
      x_1)!} = \frac{M}{x_1!(N - x_1)!}c_1^{x_1}(1 - c_1)^{N - x_1},
  \end{equation*}
  for $0\le x_1\le N$.  After setting $M = N!$, we see that $X_1$ is
  binomially distributed.  Similarly,
  \begin{equation*}
    \pi_2(x_2) = \binom{N}{x_2} c_2^{x_2}(1 - c_2)^{N-x_2},
  \end{equation*}
  for $0\le x_2\le N$.
  Therefore, we trivially have that $P(X_1 = N) = c_1^N$ and $P(X_2 =
  N) = c_2^N$, but $P(X_1 = N, X_2 = N) = 0 \ne c_1^N c_2^N$, and so
  $X_1$ and $X_2$ are not independent.
  \label{ex:no_indep}
\end{example}
\begin{remark}
  The conclusion of the previous example, that independence does not
  follow from the existence of a product-form stationary distribution,
  extends trivially to any network with a conservation relation among
  the species.
\end{remark}

\begin{example}\textbf{(First order reaction networks)} The results
  presented below for first order reaction networks are known in both
  the queueing theory and mathematical chemistry literature.  See, for
  example, \cite{Kelly1979} and \cite{Othmer2005}.  We present them
  here to point out how they follow directly from Theorem
  \ref{thm:prodform_def0}.

  We begin by defining $|v| = \sum_i v_i$ for any vector $v \in
  \R^m_{\ge 0}$.  We say a reaction network is a {\em first order
    reaction network} if $|\nu_k| \in \{0,1\}$ for each complex $\nu_k
  \in \mathcal{C}$.  Therefore, a network is first order if each entry
  of the $\nu_k$ are zeros or ones, and at most one entry can be a
  one.  It is simple to show that first order reaction networks
  necessarily have a deficiency of zero.  Therefore, the results of
  this paper are applicable to all first order reaction networks that
  are weakly reversible.  Consider such a reaction network with only
  one linkage class (for if there is more than one linkage class we
  may consider the different linkage classes as distinct networks).
  We say that the network is {\em open} if there is at least one
  reaction, $\nu_k \to \nu_k'$, for which $\displaystyle \nu_k = \vec
  0$.  Hence, by weak reversibility, there must also be a reaction for
  which $\nu_k' = \vec 0$.  If no such reaction exists, we say the
  network is {\em closed}.  If the network is open we see that $S =
  \mathbb{R}^m$, $\Gamma = \mathbb{Z}^m_{\ge 0}$ is irreducible, and
  so by Theorem \ref{thm:prodform_def0} the unique stationary
  distribution is
  \begin{equation*}
    \pi(x) = \prod_{i = 1}^m \frac{c_i^{x_i}}{x_i!} e^{-c_i}, \qquad x
    \in \Z^m_{\ge 0}, 
  \end{equation*}
  where $c \in \mathbb{R}^m_{>0}$ is the complexed balanced
  equilibrium of the associated (linear) deterministic system.
  Therefore, when in distributional equilibrium, the species numbers
  are independent and have Poisson distributions.  Note that neither
  the independence nor the Poisson distribution resulted from the fact
  that the system under consideration was a first order system.
  Instead both facts followed from $\Gamma$ being all of
  $\mathbb{Z}^m_{\ge 0}$.

  In the case of a closed, weakly reversible, single linkage class,
  first order reaction network, it is easy to see that there is a
  unique conservation relation $X_1(t) + \cdots + X_m(t) = N$, for
  some $N$.  Thus, in distributional equilibrium $X(t)$ has a
  multinomial distribution.  That is for any $x \in \Z^m_{\ge 0}$
  satisfying $x_1 + x_2 + \cdots + x_m = N$
  \begin{equation}
    \pi(x) = \binom{N}{x_1,x_2,\dots,x_m}c^{x} = \frac{N!}{x_1!\cdots
      x_m!}c_1^{x_1}\cdots c_m^{x_m},
    \label{eq:closed_dist}
  \end{equation}
  where $c \in \mathbb{R}^m_{>0}$ is the equilibrium of the associated
  deterministic system normalized so that $\sum_i c_i = 1$.
  As in the case of the open network, we note that the form of the
  equilibrium distribution does not follow from the fact that the
  network only has first order reactions.  Instead
  \eqref{eq:closed_dist} follows from the structure of the closed,
  irreducible communicating classes.
  \label{example:first_order}
\end{example}

\begin{example}\textbf{(Enzyme kinetics I)}
  Consider the possible model of enzyme kinetics given by
  \begin{equation}
    E + S \ \rightleftarrows \ ES \ \rightleftarrows \ E + P ,
    \qquad   E \ \rightleftarrows  \ \emptyset \ \rightleftarrows S ,
    \label{network:enzyme1}
  \end{equation}
  where $E$ represents an enzyme, $S$ represents a substrate, $ES$
  represents an enzyme-substrate complex, $P$ represents a product,
  and some choice of rate constants has been made.  We note that both
  $E$ and $S$ are being allowed to enter and leave the system.

  The network \eqref{network:enzyme1} is reversible and has six
  complexes and two linkage classes.  The dimension of the
  stoichiometric subspace is readily checked to be four, and so the
  network has a deficiency of zero. Theorem \ref{thm:prodform_def0}
  applies and so the stochastically modeled system has a product-form
  stationary distribution of the form \eqref{eq:stationary_product}.
  Ordering the species as $X_1 = E$, $X_2 = S$, $X_3 = ES$, and $X_4 =
  P$, the reaction vectors for this system include
  \begin{equation*}
    \left\{ \left[\begin{array}{c}
        1\\
        0\\
        0\\
        0\\
      \end{array} \right], \
    \left[\begin{array}{c}
        0\\
        1\\
        0\\
        0\\
      \end{array} \right], \
    \left[\begin{array}{c}
        -1\\
        -1\\
        1\\
        0\\
      \end{array} \right], \
    \left[\begin{array}{c}
        1\\
        0\\
        -1\\
        1\\
      \end{array} \right] \right\}.
  \end{equation*}
  We therefore see that $\Gamma = \Z^4_{\ge 0}$ is the unique closed,
  irreducible communicating class of the stochastically modeled system
  and Theorem \ref{thm:prodform_def0} tells us that in distributional
  equilibrium the species numbers are independent and have Poisson
  distributions with parameters $c_i$, which are the
  complex balanced equilibrium values of the analogous deterministically modeled
  system.
  \label{ex:first_enzyme}
\end{example}

\begin{example}\textbf{(Enzyme kinetics II)}
  Consider the possible model for enzyme kinetics given by
  \begin{equation}
    E + S \ \underset{k_{-1}}{\overset{k_1}{\rightleftarrows}} \ \ ES
    \ \ \underset{k_{-2}}{\overset{k_2}{\rightleftarrows}} \ \ E + P
    , \qquad  \emptyset  \ \
    \underset{k_{-3}}{\overset{k_3}{\rightleftarrows}} \ \
    E, 
    \label{network:enzyme2}
  \end{equation}
  where the species $E,$ $S$, $ES$, and $P$ are as in Example
  \ref{ex:first_enzyme}.  We are now allowing only the enzyme $E$ to
  enter and leave the system.  The network is reversible, there are
  five complexes, two linkage classes, and the dimension of the
  stoichiometric compatibility class is three.  Therefore, Theorem
  \ref{thm:prodform_def0} implies that the stochastically modeled
  system has a product-form stationary distribution of the form
  \eqref{eq:stationary_product}.  The only conserved quantity of the
  system is $S + ES + P$, and so $X_2(t) + X_3(t) + X_4(t) = N$ for
  some $N >0$ and all $t$.  Therefore, after solving for the
  normalizing constant, we have that for any $x \in \Z^4_{\ge 0}$
  satisfying $x_2 + x_3 + x_4 = N$
  \begin{equation*}
    \pi(x) = e^{-c_1}\frac{c_1^{x_1}}{x_1!} \frac{N!}{x_2!x_3!x_4!}
    c_2^{x_2} c_3^{x_3} c_4^{x_4} = e^{-c_1}\frac{c_1^{x_1}}{x_1!}
    \binom{N}{x_2,x_3,x_4} c_2^{x_2} c_3^{x_3} c_4^{x_4},
  \end{equation*}
  where $c = (k_3/k_{-3},c_2,c_3,c_4)$ has been chosen so that $c_2 +
  c_3 + c_4 = 1$.  Thus, when the stochastically modeled system is in
  distributional equilibrium we have that: (a) $E$ has a Poisson
  distribution with parameter $k_3/k_{-3}$, (b) $S, \ ES,$ and $P$ are
  multinomially distributed, and (c) $E$ is independent from $S, \
  ES,$ and $P$.
  \label{ex:second_enzyme}
\end{example}

\subsection{The multiscale nature of reaction networks}
\label{sec:multiscale}

Within a cell, some chemical species may be present in much greater
abundance than others. In addition, the rate constants $\kappa_k$ may
vary over several orders of magnitude. Consequently, the scaling limit
that gives the classical deterministic law of mass action detailed in Section \ref{sec:scaling} may not be appropriate, and a different approach to deriving a scaling limit
approximation for the basic Markov chain model must be considered.
As a consequence of the multiple scales in a network model, it may
be possible to separate the network into subnetworks of species and
reactions, each dominated by a time scale of a specific
magnitude. Within each subnetwork, the graph structure and
stoichiometric properties may determine properties of the asymptotic
solutions of the subnetwork.

\begin{example}

Consider the reaction network 
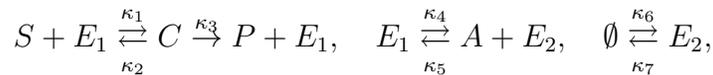
\begin{equation*}
  S + E_1 \overset{\kappa_1}{\underset{\kappa_2}{\rightleftarrows}} C \overset{\kappa_3}{\rightarrow} P + E_1, \quad E_1
  \overset{\kappa_4}{\underset{\kappa_5}{\rightleftarrows}} A + E_2, \quad \emptyset \overset{\kappa_6}{\underset{\kappa_7}{\rightleftarrows}} E_2,
\end{equation*}
where $\emptyset \to E_2$ and $E_2 \to \emptyset$ represent production
and degradation of $E_2$, respectively, $S$ is a substrate being
converted to a product $P$, $E_1$ and $E_2$ are enzymes, and $A$ is a
substrate that reacts with $E_2$ allosterically to transform it into
an active form.

We suppose that $(i)$ the enzymes $E_1$, $E_2$ and the substrate $A$ are in
relatively low abundances, $(ii)$ the substrate $S$ has a large
abundance of ${\cal O}(V)$, and $(iii)$ the reaction rates are also of
the order ${\cal O}(V)$.  We change notation slightly and denote the number of molecules of species $A$ at time $t$ as $X^V_A(t)$, and similarly for the other species.  Further, we denote $X^V_S(t)/V = Z^V_S(t)$. Combined with the conservation relation $X_{E_1}^V + X^V_{C} + X^V_A = M \in \Z_{> 0}$, the scaled equations for the stochastic model are
\begin{align*}
  &Z^V_S(t) = Z^V_S(0) - V^{-1}
  Y_1(V\int_0^t\kappa_1Z^V_S(s)X^V_{E_1}(s)ds) + V^{-1} Y_2(V\int_0^t
  \kappa_2 X^V_{C}(s)ds) \label{eq:slow}\\
  &X^V_{E_1}(t) =
  X^V_{E_1}(0)-Y_1(V\int_0^t\kappa_1Z^V_S(s)X^V_{E_1}(s)ds) +
  Y_2(V \int_0^t\kappa_2X^V_{C}(s)ds)\notag \\
  & \hspace{.6in}+ Y_3(V\int_0^t\kappa_3X^V_{C}(s)ds) -
  Y_4(V\int_0^t\kappa_4X^V_{E_1}(s) ds ) + Y_5(V\int_0^t\kappa_5
  X^V_A(s)X^V_{E_2}(s)ds)\notag \\
  &X^V_A(t) = X_A^V(0) + Y_4(V\int_0^t\kappa_4X^V_{E_1}(s)ds) -
  Y_5(V\int_0^t \kappa_5X^V_A(s) X_{E_2}^V(s) ds )\notag \\
  &X_{E_2}^V(t) = X_{E_2}^V(0) + Y_6(V\kappa_6t)+Y_4(V \int_0^t
  \kappa_4 X^V_{E_1}(s) ds ) - Y_5
  (V\int_0^t\kappa_5X^V_A(s)X_{E_2}^V(s) ds)\notag \\
  &\hspace{.6in} -Y_7 (V\int_0^t\kappa_7X^V_{E_2}(s)ds),\notag
\end{align*}
where the $Y_i$ are unit-rate Poisson processes.
The first equation satisfies
\begin{equation*}
  Z^V_S(t) = Z^V_S(0) - V^{-1}
  Y_1(V\int_0^t  \kappa_1Z^V_S(s)  \int_{-\infty}^{\infty} x
  \mu^V_{s}(dx) ds) + V^{-1} Y_2(V\int_0^t \kappa_2
  \int_{-\infty}^{\infty} x \eta^V_{s}(dx) ds), 
\end{equation*}
where $\mu^V_s(A) = I_{\{X_{E_1}^V(s) \in A\}}$ and $\eta^V_s(A) =
I_{\{X_{C}^V(s) \in A\}}$ are the respective occupation measures.
Using methods from stochastic averaging (see, for example,
\cite{Ball06, Kurtz92book}), as $V \to \infty$ the fast system is
``averaged out:''
\begin{equation}
  Z_S(t) = Z_S(0) -  \int_0^t  \kappa_1Z_S(s)  \int_{-\infty}^{\infty}
  x \mu_{s}(dx) ds +  \int_0^t \kappa_2 \int_{-\infty}^{\infty} x
  \eta_{s}(dx) ds, 
  \label{eq:need_mean}
\end{equation}  
where $\mu_s$ and $\eta_s$ are the stationary distributions of
$X_{E_1}$ and $X_{C}$, respectively, of the fast subsystem with
$Z_S(s)$ held constant (assuming a stationary distribution
exists). This reduced network 
(i.e. the fast subsystem) is
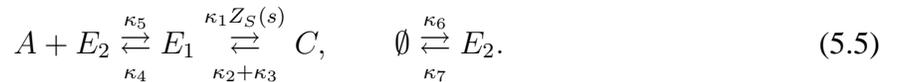
\begin{equation}
  A + E_2 \overset{\kappa_5}{\underset{\kappa_4}{\rightleftarrows}}
  E_1 \overset{\kappa_1Z_S(s)}{\underset{\kappa_2 +
      \kappa_3}{\rightleftarrows}} C, \qquad \emptyset
  \overset{\kappa_6}{\underset{\kappa_7}{\rightleftarrows}} E_2.
  \label{eq:fast}
\end{equation}
Setting $z = Z_S(s)$ we have the following equilibrium relations for
the moments of the above network
\begin{align}
\begin{split}
  \kappa_4\E[X_{E_1}] - \kappa_5\E[X_A X_{E_2}] &= 0\\
  -(\kappa_1z+\kappa_4)\E[X_{E_1}]+(\kappa_2+\kappa_3)\E[X_{C}]+\kappa_5
  \E[X_{A}X_{E_2}]& = 0\\
  \kappa_6+\kappa_4\E[X_{E_1}]-\kappa_5\E[X_{A}X_{E_2}]-\kappa_7\E[X_{E_2}]&=0\\
  \E[X_{E_1}] + \E[X_{C}]+ \E[X_A]&=M.
  \end{split}
  \label{eq:solving}
\end{align}
$\E[X_{E_1}]$ and $\E[X_C]$, which are both functions of $z$ and
needed in equation \eqref{eq:need_mean}, can not be explicitly solved
for via the above equations without extra tools as
\eqref{eq:solving} is a system of four equations with five unknowns.
This situation arises frequently as it stems from the nonlinearity of
the system.  However, the
network \eqref{eq:fast} consists of five complexes, two connected
components, and the dimension of its stoichiometric subspace is three.
Therefore, its deficiency is zero.  As it is clearly weakly
reversible, Theorem \ref{thm:prodform_main} applies and, due to the
product form of the distribution and the unboundedness of the support
of $X_{E_2}$, it is easy to argue that $X_{E_2}$ is independent of
$X_A$, $X_{E_1}$, and $X_{C}$ when in equilibrium.  Thus,
$\E[X_AX_{E_2}] = \E[X_A]\E[X_{E_2}]$ and the first moments can be
solved for as functions of $Z_S(s)$.  After solving and inserting
these moments, \eqref{eq:need_mean} becomes
\begin{align*}
  Z_S(t) &= Z_S(0) - \int_0^t \frac{\kappa_1 \kappa_3 \kappa_5
    \kappa_6 M Z_S(s)}{(\kappa_5\kappa_6 + \kappa_7\kappa_4)(\kappa_2
    + \kappa_3) + \kappa_1 \kappa_5 \kappa_6 Z_S(s)} ds,
\end{align*}
which is Michaelis-Menten kinetics.  
\label{ex:multiscale}
\end{example}

\setcounter{equation}{0}
\section{More general kinetics}
\label{sec:result_Gen}

In this section we extend our results to systems with more general
kinetics than stochastic mass action.  The generalizations we make are
more or less standard for the types of results presented in this paper
(see, for example, \cite{Kelly1979}, Section 8.5, \cite{Whittle86},
Chapter 9).  What is surprising, however, is that the conditions of
the deficiency zero theorem of Feinberg (which are conditions on
mass-action deterministic systems) are also sufficient to guarantee
the existence of stationary distributions of stochastically modeled
systems even when the intensity functions are not given by
\eqref{eq:stoch_MA}.  It is interesting to note that the
generalizations made here for the stochastic deficiency zero Theorem
\ref{thm:prodform_def0} are similar to those made in
\cite{Sontag2001}, which generalized Feinberg's deficiency zero
Theorem \ref{thm:def0} in the deterministic setting.

\vspace{.125in}

Suppose that the intensity functions of a stochastically modeled
system are given by
\begin{equation}
  \lambda_k(x) = \kappa_k \prod_{i = 1}^m \prod_{j=0}^{\nu_{ik}-1}
  \theta_i(x_i - j) = \kappa_k \prod_{i = 1}^m 
  \theta_i(x_i)\theta_i(x_i - 1)\theta_i(x_i - (\nu_{ik} - 1)),
  \label{eq:gen_kin}
\end{equation}
where the $\kappa_k$ are positive constants, $\theta_i : \Z \to
\R_{\ge 0}$, $\theta_i(x) = 0$ if $x \le 0$, and we use the convention
that $\prod_{j=0}^{-1}a_j = 1$ for any $\{a_j\}$.  Note that the final
condition allows us to drop the indicator functions of
\eqref{eq:stoch_MA}.  As pointed out in \cite{Kelly1979}, the function
$\theta_i$ should be thought of as the ``rate of association'' of the
$i$th species.  We give a few interesting choices for $\theta_i$.  If
$\theta_i(x_i) = x_i$ for $x_i \ge 0$, then \eqref{eq:gen_kin} is
stochastic mass-action kinetics.  However, if for $x_i \ge 0$
\begin{equation}
  \theta_i(x_i) = \frac{v_ix_i}{k_i + x_i},
  \label{eq:MM}
\end{equation}
for some positive constants $k_i$ and $v_i$, then the system has a
type of stochastic Michaelis-Menten kinetics (\cite{keener1998},
Chapter 1).  Finally, if $|\nu_k| \in \{0,1\}$ and $\theta_i(x_i) =
\min\{n_i,x_i\}$ for $x_i \ge 0$, then the dynamical system models an
$M/M/n$ queueing network in which the $i$th species (and in this case
complex) represents the queue length of the $i$th queue, which has
$n_i$ servers who work on a first come, first serve basis.

The main restriction imposed by \eqref{eq:gen_kin} is that for any
reaction for which the $i$th species appears in the source complex,
the rate of that reaction must depend upon $X_i$ via $\theta_i(X_i)$
only.  Therefore, if, say, the $i$th species is governed by the
kinetics \eqref{eq:MM}, then the constants $k_i$ and $v_i$ must be the
same for each intensity which depends upon $X_i$ (although the $v_i$
may be incorporated into the rate constants $\kappa_k$, and so the
real restriction is on the constant $k_i$).  However, systems with
intensities given by \eqref{eq:gen_kin} are quite general in that
different kinetics can be incorporated into the same model through the
functions $\theta_i$.  For example, if in a certain system species
$S_1$ is modeled to be governed by Michaelis-Menten kinetics
\eqref{eq:MM} and species $S_2$ is modeled to be governed by
mass-action kinetics, then the reaction $S_1 + S_2 \to \nu_k'$ would
have intensity
\begin{equation*}
  \lambda_k(x) = \kappa_k \frac{v_1x_1}{k_1 + x_1}x_2,
\end{equation*}
for some constant $\kappa_k$.

In following we use the convention that $\prod_{j=1}^0 a_j = 1$ for
any choice of $\{a_j\}$.

\begin{theorem}
  Let $\{ \S, \C, \Re \}$ be a stochastically modeled chemical
  reaction network with intensity functions \eqref{eq:gen_kin}.
  Suppose that the associated mass-action deterministic system with
  rate constants $\{\kappa_k\}$ has a complex balanced equilibrium $c \in
  \R^m_{>0}$.  Then the stochastically modeled system admits the stationary
  distribution 
  \begin{equation}
  	\pi(x) = M\prod_{i = 1}^m
        \frac{c_i^{x_i}}{\prod_{j=1}^{x_i}\theta_i(j)}, \qquad x \in
        \Z^m_{\ge 0}, 
	\label{eq:pfrm_gen}
  \end{equation}
  where $M>0$ is a normalizing constant, provided that
  \eqref{eq:pfrm_gen} is summable.  If $\Z^m_{\ge 0}$ is irreducible,
  then \eqref{eq:pfrm_gen} is the unique stationary distribution,
  whereas if $\Z^m_{\ge 0}$ is not irreducible then the $\pi_{\Gamma}$
  of equation \eqref{eq:pi_gamma} are given by the product-form
  stationary distributions
  \begin{equation}
    \pi_{\Gamma}(x) = M_{\Gamma} \prod_{i = 1}^m
    \frac{c_i^{x_i}}{\prod_{j=1}^{x_i}\theta_i(j)}, \qquad x \in
    \Gamma, 
    \label{eq:pfrm_gen2}
  \end{equation}
  and $\pi_{\Gamma}(x) = 0$ otherwise, where $M_{\Gamma}>0$ is a
  normalizing constant, provided that \eqref{eq:pfrm_gen2} is
  summable.
  \label{thm:firstgen}
\end{theorem}

\begin{proof}
  The proof consists of plugging \eqref{eq:pfrm_gen} and
  \eqref{eq:gen_kin} into equation \eqref{eq:stationary2} and
  verifying that $c$ being a complex balanced equilibrium is
  sufficient.  The details are similar to before and so are omitted.
\end{proof}

\begin{remark}
  We simply remark that just as Theorem \ref{thm:prodform_def0}
  followed directly from Theorem \ref{thm:prodform_main}, the results
  of Theorem \ref{thm:firstgen} hold, independent of the choice of
  rate constants $\kappa_k$, so long as the associated network is
  weakly reversible and has a deficiency of zero.
  \label{rmk:gen_def0}
\end{remark}

\begin{example}
  Consider a network, $\{\S,\C,\Re\}$, that is weakly reversible and
  has a deficiency of zero.  Suppose we have modeled the dynamics
  stochastically with intensity functions given by \eqref{eq:gen_kin}
  with each $\theta_i$ given via \eqref{eq:MM} for some choice of
  $v_i>0$ and $k_i$ a nonnegative integer.  That is, we consider a
  system endowed with stochastic Michaelis-Menten kinetics. Then,
  \begin{equation*}
    \prod_{j=1}^{x_i}\theta_i(j) = \prod_{j=1}^{x_i} \frac{v_i j}{k_i
      + j} = v_i^{x_i}/\binom{k_i +x_i}{x_i}. 
  \end{equation*}
  Thus, our candidate for a stationary distribution is
  \begin{equation}
    \pi (x) = M \prod_{i = 1}^m
    \frac{c_i^{x_i}}{\prod_{j=1}^{x_i}\theta_i(j)} = M \prod_{i = 1}^m
    \binom{k_i + x_i}{x_i}\left(\frac{c_i}{v_i}\right)^{x_i}. 
    \label{eq:stat_MM}
  \end{equation}
  Noting that
  \begin{equation*}
    \binom{k_i + x_i}{x_i} = O(x_i^{k_i}), \hspace{.1in} x_i \to \infty,
  \end{equation*}
  we see that $\pi(x)$ given by \eqref{eq:stat_MM} is summable if $c_i
  < v_i$ for each species $S_i$ whose possible abundances are
  unbounded.  In this case, \eqref{eq:stat_MM} is indeed a stationary
  distribution for the system.  We note that the condition $c_i < v_i$
  for each species $S_i$ is both necessary and sufficient to guarantee
  summability if $Z^m_{\ge 0}$ is irreducible, as in such a situation the
  species numbers are independent.
  \label{ex:MM}
\end{example}

\begin{example}
  In \cite{LevineHwa2007}, Levine and Hwa computed and analyzed the
  stationary distributions of different stochastically modeled
  chemical reaction systems with Michaelis-Menten kinetics
  \eqref{eq:MM}.  The models they considered included among others:
  directed pathways ($\emptyset \to S_1 \to S_2 \to \cdots \to S_L \to
  \emptyset$), reversible pathways ($\emptyset \to S_1
  \leftrightarrows S_2 \leftrightarrows \cdots \leftrightarrows S_L
  \to \emptyset$), pathways with dilution of intermediates ($S_i \to
  \emptyset$), and cyclic pathways ($S_L \to S_1$).  Each of the
  models considered in \cite{LevineHwa2007} is biologically motivated
  and has a first order reaction network ($|\nu_k| \in \{0,1\}$, see
  Example \ref{example:first_order}), which guarantees that they have
  a deficiency of zero.  Further, the networks of the models
  considered are weakly reversible; therefore, the results of the
  current paper, and in particular Theorem \ref{thm:firstgen} and the
  remark that follows, apply so long as the restrictions discussed in
  the paragraph preceding Theorem \ref{thm:firstgen} are met.  While
  these restriction are not always met (for example, dilution is
  typically modeled with a linear intensity function and there is no
  reason for the $k_i$ of a forward and a backward reaction for a
  species $S_i$ in a reversible pathway to be the same), they found
  that the stationary distributions for these models are either of
  product form (when the restrictions are met) or near product form
  (when the restrictions are not met).  Further, because $\Z^m_{\ge
    0}$ is irreducible in each of these models, the product form of
  the distribution implies that the species numbers are independent.
  It is then postulated that the independence of the species numbers
  could play an important, beneficial, biological role (see
  \cite{LevineHwa2007} for details).  Similar to the conclusions we
  drew in Example \ref{example:first_order}, Theorem
  \ref{thm:firstgen} and the remark that follows point out how the
  models analyzed in \cite{LevineHwa2007} are actually special cases
  of a quite general family of systems that have both the product form
  and independence properties, and that these properties may be more
  widespread, and taken advantage of by living organisms, than
  previously thought.
\end{example}

We return to the result of Example \ref{ex:MM} pertaining to the
summability of \eqref{eq:stat_MM} and show that this can be
generalized in the following manner.

\begin{theorem}
  Suppose that for some closed, irreducible $\Gamma \subset \Z^m_{\ge
    0}$, $\pi_{\Gamma}: \Gamma \to \R_{\ge 0}$ satisfies
  \begin{equation*}
    \pi_{\Gamma}(x) = M\prod_{i = 1}^m
    \frac{c_i^{x_i}}{\prod_{j=1}^{x_i}\theta_i(j)}, 
  \end{equation*}
  for some $c \in \R^m_{> 0}$ and $M > 0$, where $\theta_i: \Z_{\ge
    0} \to \R_{\ge 0}$ for each $i$.  Then $\pi_{\Gamma}(x)$ is
  summable if for each $i$ for which $\sup \{x_i \ | \ x \in \Gamma\}
  = \infty$ we have that $\theta_i(j) > c_i + \epsilon$ for some
  $\epsilon > 0$ and $j$ sufficiently large.  \label{thm:summability}
\end{theorem}

\begin{proof}
  The conditions of the theorem immediately imply that there are
  positive constants $C$ and $\rho$ for which $\pi_{\Gamma}(x) <
  Ce^{-\rho|x|}$, for all $x \in \Gamma$, which implies that
  $\pi_{\Gamma}(x)$ is summable.
\end{proof}

It is tempting to believe that the conditions of Theorem
\ref{thm:summability} are in fact necessary, as in the case when
$Z^m_{\ge 0}$ is irreducible.  The following simple example shows this
not to be the case.

\begin{example}
  Consider the reaction system with network 
  \begin{equation*}
    \emptyset  \rightleftarrows  S_1 + S_2,
  \end{equation*}
  where the rate of the reaction $\emptyset \to S_1 + S_2$ is
  $\lambda_1(x) = 1$, and the rate of the reaction $S_1 + S_2 \to
  \emptyset$ is $\lambda_2(x) = 1\times \theta_1(x_1)\theta_2(x_2)$,
  where
  \begin{equation*}
    \theta_1(x_1) = \frac{3x_1}{1 +x_1}, \quad \theta_2(x_2) =
    \frac{(1/2)x_2}{1 +x_2}. 
  \end{equation*}
  Assume further that $X_1(0) = X_2(0)$.  For the more physically
  minded readers, we note that this model could describe a reaction
  system for which there is a chemical complex $C = S_1S_2$ that
  sporadically breaks into its chemical constituents, which may then
  re-form.  The complex $C$ may be present in such high numbers
  relative to free $S_1$ and $S_2$ that we choose to model it as
  fixed, which leads to the above reaction network.
  
  We note that in this case, the reaction rates $\{\kappa_k\}$ for the
  corresponding mass-action deterministic system are both equal to
  one, and so an equilibrium value guaranteed to exist for the
  deterministically modeled system by the deficiency zero theorem is
  $c = (1,1)$.  This system does not satisfy the assumptions of
  Theorem \ref{thm:summability} because both $X_1$ and $X_2$ are
  unbounded and $\lim_{j \to \infty}\theta_2(j) = 1/2 < 1 = c_2$.
  However, for any $x \in \Gamma = \{x \in \Z^2_{\ge 0} \ : \ x_1 =
  x_2\},$
  \begin{equation*}
    \pi_{\Gamma}(x) = \binom{1 +
      x_1}{x_1}\left(\frac{1}{3}\right)^{x_1} \binom{1 +
      x_2}{x_2}\left(\frac{1}{(1/2)}\right)^{x_2} = \binom{1 +
      x_1}{x_1}^2\left(\frac{2}{3}\right)^{x_1},   
  \end{equation*}
  which is summable over $\Gamma$.
  \label{ex:not_necessary}
\end{example}

For the most general kinetics handled in this paper, we let the
intensity functions of a stochastically modeled system be given by
\begin{equation}
  \lambda_k(x) = \kappa_k \frac{\theta(x)}{\theta(x - \nu_k)}
  \prod_{\ell = 1}^m 1_{\{x_{\ell} \ge \nu_{\ell k}\}}, 
  \label{eq:most_gen_kin}
\end{equation}
where the $\kappa_k$ are positive constants, and $\theta: \Z^m \to
\R_{> 0}$.  Note that if
\begin{equation*}
  \theta(x) = \prod_{i = 1}^m \prod_{j=1}^{x_i}\theta_i(j),
\end{equation*}
for some functions $\theta_i$, then \eqref{eq:most_gen_kin} is
equivalent to \eqref{eq:gen_kin}, and so the following theorem implies
Theorem \ref{thm:firstgen}.  It's proof is similar to the previous
theorems and so is omitted.

\begin{theorem} 
  Let $\{ \S, \C, \Re \}$ be a stochastically modeled chemical
  reaction network with intensity functions \eqref{eq:most_gen_kin}.
  Suppose that the associated mass-action deterministic system with
  rate constants $\{\kappa_k\}$ has a complex balanced equilibrium $c \in
  \R^m_{>0}$.  Then the stochastically modeled system admits the stationary
  distribution 
  \begin{equation}
  	\pi(x) =  M\frac{1}{\theta(x)}\prod_{i = 1}^m
    c_i^{x_i}, \quad x \in \Z^m_{\ge 0},
    \label{eq:pfrm_most_gen1}
  \end{equation}
  where $M>0$ is a normalizing constant, provided that
  \eqref{eq:pfrm_most_gen1} is summable.  If $\Z^m_{\ge 0}$ is
  irreducible, then \eqref{eq:pfrm_most_gen1} is the unique stationary
  distribution, whereas if $\Z^m_{\ge 0}$ is not irreducible then the
  $\pi_{\Gamma}$ of equation \eqref{eq:pi_gamma} are given by the
  product-form stationary distributions
  \begin{equation}
    \pi_{\Gamma}(x) = M_{\Gamma} \frac{1}{\theta(x)}\prod_{i = 1}^m
    c_i^{x_i}, \qquad x \in \Gamma,
    \label{eq:pfrm_most_gen2}
  \end{equation}
  and $\pi_{\Gamma}(x) = 0$ otherwise, where $M_{\Gamma}>0$ is a
  normalizing constant, provided that \eqref{eq:pfrm_most_gen2} is
  summable.
  \label{thm:thirdgen}
\end{theorem}

\begin{remark}
  Similar to the remark following Theorem \ref{thm:firstgen}, we point
  out that the results of Theorem \ref{thm:thirdgen} hold, independent
  of the choice of rate constants $\kappa_k$, so long as the
  associated network is weakly reversible and has a deficiency of
  zero.
  \label{rmk:_mostgen_def0}
\end{remark}

\vspace{.25in}
\begin{center}
  {\bf \large Acknowledgments}
\end{center}

We gratefully acknowledge the financial support of the National
Science Foundation through grant NSF-DMS-0553687.

\bibliographystyle{amsplain} \bibliography{prodform}

\end{document}